\documentclass[preprint,12pt]{elsarticle}

\makeatletter
\def\ps@pprintTitle{%
 \let\@oddhead\@empty
 \let\@evenhead\@empty
 \def\@oddfoot{}%
 \let\@evenfoot\@oddfoot}
\makeatother

\usepackage{amssymb}
\usepackage{amsmath}%
\newdefinition{definition}{Definition}
\newtheorem{proposition}{Proposition}
\newtheorem{theorem}{Theorem}
\newtheorem{example}{Example}
\newtheorem{lemma}{Lemma}
\newproof{pf}{Proof}
\usepackage{booktabs,tabu}
\begin{document}

\begin{frontmatter}



\title{Constraints and Conditions: the Lasso Oracle-inequalities}


%\author{Niharika Gauraha \\ Indian Statistical Institute}

%\address{}
 %\author[label1,label2]{N and S}
 %\address[label1]{Indian Statistical Institute}
 %\address[label2]{Indian Statistical Institute}
 
\author[add1]{Niharika Gauraha}
\ead{niharika.gauraha@gmail.com}
\author[add2]{Swapan K. Parui}
\ead{swapan.parui@gmail.com}
\address[add1]{Systems Science and Informatics Unit\\
       Indian Statistical Institute\\
       8th Mile, Mysore Road Bangalore, India}
\address[add2]{Computer Vision and Pattern Recognition Unit\\
       Indian Statistical Institute\\
       203 B.T. Road Kolkata, India\\}

\begin{abstract}
We study various constraints and conditions on the true coefficient vector and on the design matrix to establish non-asymptotic oracle inequalities for the prediction error, estimation accuracy and variable selection for the Lasso estimator in high dimensional sparse regression models. We review results from the literature and we provide simpler and detailed derivation for several boundedness theorems. In addition, we complement the theory with illustrated examples. %We review results from the literature and we provide simpler proofs and lower bounds for several boundedness theorems. In addition, we complement the theory with illustrated examples.

\end{abstract}


\end{frontmatter}


\section{Introduction}
We consider the usual linear regression model

\begin{align} \label{eq:lr} 
	\textbf{Y} &=  \textbf{X} \beta^{0}+\epsilon,
\end{align}
with response vector $Y_{n \times 1}$, design matrix $X_{ n \times p}$, true underlying coefficient vector $\beta^0_{p \times 1}$ and error vector $\epsilon_{n\times 1}$.  Mainly, we focus on the models where the number of predictors (p) is much larger than the number of observations (n), $p \gg n$, and we assume sparsity in  $\beta^{0}$. A lot of research has been devoted  to penalized estimators over a more than decade, we refer to a few papers here \cite{Tibshirani}, \cite{Zou}, \cite{hui}, \cite{Candes}, \cite{Nicolai2} and  \cite{Belloni} etc. The Least Absolute Shrinkage and selection operator (Lasso) is a popular penalized method for simultaneous estimation and variable selection. We consider the Lasso for high dimensional sparse regression models and we study various conditions on the design matrix required to establish non-asymptotic oracle inequalities. 
%We study some sufficient conditions on the design matrix required to establish oracle inequalities for the Lasso in regression. We usually consider the Gaussian regression model with the fixed design matrix. First we derive the oracle inequalities for the noiseless case.

The oracle property means that the penalized estimator is
asymptotically equivalent to the oracle estimator that is as good as the true underlying model was given in advance, i.e. $\beta^0 = \beta^{0}_{S_0}$, where $S_0$ is the true support set. Oracle results for estimation and prediction has been established assuming restricted eigenvalue conditions \cite{Bickel} and its various forms, see \cite{Sara2}, \cite{Zhang}, \cite{Buhlmann1}, \cite{Sara1}, \cite{Sara3} and \cite{Nicolai}. For variable selection consistency, the design matrix $\textbf{X}$ must satisfy irrepresentable condition(\cite{Zhao}) or neighborhood stability condition(\cite{Buhlmann4}) along with beta-min conditions.

Our aim is two fold: to review various results from the literature, and to provide simpler and detailed derivation for several boundedness theorems. We also give illustrated examples to support the theory.

The rest of this paper is organized as follows. In Section 2, we state notations and assumptions. In section 3, we start with the sparse recovery in the simplest model where observations are noiseless , then we study sufficient conditions required on the design matrix for the exact recovery in a noiseless situation. In section 4, we provide theory with illustrated examples on the various conditions and constraints required for the Lasso oracle-inequalities to hold. In Section 5, we (re)derive some useful bounds which we require to prove bounds on various loss functions . In section 6, we derive oracle inequalities for prediction accuracy and estimation error for the Lasso for both noiseless and noisy case. We also (re)prove that the irrepresentable condition is almost necessary and sufficient for variable selection for the noiseless Lasso problem. We shall provide conclusion in section 7.

%We also show that the variable selection is consistent for high dimensional sparse problems 

\section{Notations and Assumptions}
In this section, we state notations and assumptions, and we also define required concepts, which are applied throughout this paper.

We consider the usual linear regression model as given in Eq. (\ref{eq:lr}), and we assume that the random errors are i.i.d. with $N(0, \sigma^2)$, that is  $\epsilon \sim  N_n(0, \sigma^2 I)$, the design matrix $\textbf{X}$ is fixed and $p \gg n$.\\
\\
$\ell_0$-norm counts the number of non zero coefficients, and is defined as:
\begin{align} \label{eq:l0}
	\|\beta\|_0 = \textstyle \sum_{j=1}^p (\beta_j)^0, if \beta_j \neq 0
\end{align}
$\ell_1$-norm is defined as:
\begin{align} \label{eq:l1}
	\|\beta\|_1 = \textstyle \sum_{j=1}^p |\beta_j|
\end{align}
$\ell_2$-norm squared is defined as:
\begin{align} \label{eq:l2}
	\|\beta\|^{2}_2 = \textstyle \sum_{j=1}^p \beta^{2}_{j}
\end{align}
The $\ell_\infty$-norm is defined as:
\begin{align} \label{eq:linf}
	\|\beta\|_{\infty} = \textstyle max_{1 \leq i \leq n |} |\beta_j|
\end{align}
In general, $\ell_q$-norm is defined as
\begin{align*}
    {\displaystyle \|\beta \|_{q}=\left(\sum _{j=1}^{p}|\beta_{j}|^{q}\right)^{1/q}} 
\end{align*}
%We denote the true regression vector as $\beta^0$ and 
The true active set denoted by $S_0$ is the support of the subset selection solution ($S_0 = supp(\beta_0)$) and defined as (we may also use ``S" for the fixed active set)
\begin{align} \label{eq:s0}
	S_0 &= \{j; \beta^{0}_{j} \neq 0\}.
\end{align}
We denote a Lasso estimated parameter vector as $\hat{\beta}$. Assume $\lambda > 0 $ is the regularization parameter, then the Lasso estimator is computed as:
\begin{align} \label{eq:lasso}
	\hat{\beta} := \hat{\beta}(\lambda)= \mathop{arg min}_{\beta \in \mathbb{R}^p} \{\frac{1}{n} \| \textbf{y} - \textbf{X} \beta \|_{2}^{2}+ \lambda * \|\beta\|_1 \}
\end{align}
and the Lasso estimated active set is denoted as $\hat{S}$ and defined as
\begin{align} \label{eq:s1}
	\hat{S} = \{j; \hat{\beta}_{j} \neq 0\}
\end{align}
For the noiseless case $Y = X\beta^0$, the Lasso problem for some fixed $\lambda > 0$, is defined as:
\begin{align} \label{noiselessLasso}
 \beta^{\star} := \mathop{argmin}_{\beta} \{ \| X\beta - X\beta^0 \|_{2}^{2} + \lambda \| \beta \|_1 \}
\end{align}
The sign function is defined as:
\begin{align}
	sign(x) = \left\lbrace \begin{array}{ll}
	-1 & \text{ if } x < 0  \\
	0 & \text{ if }  x = 0\\
	1 & \text{ if }  x > 0
\end{array}	 \right.
\end{align}
Null space of a matrix X is defined as \begin{align}
 null(X) = \{ \beta \in \mathbb{R}^p | X\beta = 0\}.
\end{align}
$\ell_1$ norm can be upper bounded by $\ell_2$ norm as (see \cite{Buhlmann1})
\begin{align} \label{eq:L12}
	\|\beta_{S_0}\|_{1} \leq \sqrt{s_0}\|\beta_{S_0} \|_{2}
\end{align}
The (scaled) Gram matrix(covariance matrix) is defined as $\hat{\Sigma}= \frac{X'X}{n}$. We may use $\Sigma$ (without the ``hat" superscript) for the covariance matrix of the fixed design $\textbf{X}$. 
The $\beta_S$ has zeroes outside the set $S$, as
\[
	\beta_S = \{ \beta_j I(j \in S) \}
\]
and \[\beta = \beta_S + \beta_{S^c}.\]
For the fixed active set $S$, the covariance matrix can be partitioned for the active and the redundant variables as
\begin{align} \label{eq:covMat}
	\Sigma = \left[ \begin{array}{cc}
	\Sigma_{11} = \Sigma(S) & \Sigma_{12}(S)\\
	\Sigma_{21}(S)\quad \quad &\Sigma_{22} = \Sigma(S^c)	
	\end{array} \right].
\end{align}
We also assume that $\Sigma_{11}$ is non-singular, that is $\Lambda_{min}(\Sigma_{11}) > 0$.
%squared $L_2(Q)$ norm is defined as:	\| x \|^{2}_{2} = x'x
%Signal to Noise Ratio(SNR) is defined as
%\begin{align*}
% 	SNR &= \frac{ \|X \beta^0 \|_2 } {\sqrt{n}\sigma} % %\sqrt{\frac{Var(X\beta^0)}{\sigma^2}}
% \end{align*}
The Lasso error is denotes as (for simplicity we may omit the superscript hat notation)
\[
	\hat{\Delta}= \hat{\beta} - \beta^0
\]
 Soft thresholding function is defined as follows.
 \begin{align*}
 S_{\lambda}(x) =  \left\lbrace \begin{array}{ll}
	sign(x)(|x| - \lambda) & \text{ if } |x| \geq \lambda  \\
	0 &  \text{ otherwise }
\end{array}	 \right.
 \end{align*}
Holder's inequality for two real function f and g, is given as 
\begin{align} \label{eq:holder}
	\|fg\|_1 \leq \|f\|_p\; &\|g\|_q, \quad \text{ where } p, q \in [1, \infty] \text{ and }\frac{1}{p}+\frac{1}{q} = 1.
\end{align}
%For $ p= \infty \; and \; q = 1 $, it can be interpreted as 	$| \sum f_i g_i |  \leq  |  \sum f_i|* max(g_j) $
	
%\section{Some conditions on Design Matrix and Loss Function}
\section{Sparse Recovery in Noiseless Case}
In order to get an idea of how some conditions on design matrix ensures the exact (or an optimal) solution, we start with the simplest model, where the observations are noiseless. We consider a linear equation as
\begin{align}\label{bplp}
	\textbf{Y} = \textbf{X} \beta^{0},
\end{align}
where $\textbf{Y, X}$ and $\beta^0$ are defined as in previous section. When $p>n$, this is an under-determined linear system, and there is no unique solution or infinitely many solutions exist. We need to pick one out of many solutions that satisfies certain property, for example, sparsity assumption for the true $\beta^0$. More concretely, we aim to obtain a solution to the linear equation (\ref{bplp}), that has the fewest number of non-zero entries in $\beta^0$.
%We use $\ell_1$ regularization, which basically leads to the following convex optimization problem(also known as the basis pursuit linear program).
This problem can be best described as $\ell_0$-norm constrained optimization problem and defined as
\begin{align} \label{noiselessL0}
	\mathop{minimize}_{\beta}\;  \| \beta \|_0 \text{ such that } \textbf{Y} = \textbf{X} \beta
\end{align}
But the above optimization problem is non-convex and NP-complete. Then, we consider a nearest convex problem, which is $\ell_1$-norm constrained convex optimization problem defined as follows.
\begin{align} \label{noiselessL1}
	\mathop{minimize}_{\beta} \; \| \beta \|_1 \text{ such that } \textbf{Y} = \textbf{X} \beta
\end{align}
It is also known as the Basis Pursuit Linear Program (BPLP), see \cite{Cohen}. Now, we state the assumption under which the solution of BPLP (\ref{noiselessL1}) is equivalent to the solution of $\ell_0$ non-convex problem (\ref{noiselessL1}). %Now, we state various constraints and conditions required for exact estimation of the $\beta^0$. %First we discuss that under assumption of the restricted null space property, the exact recovery is possible for BPLP.
\begin{theorem}
The BPLP estimates the true $\beta^0$ exactly if and only if the design matrix \textbf{X} satisfies the restricted null space property with respect to the true active set $S_0$.
\end{theorem}
For proof of the above theorem we refer to \cite{Cohen}. The restricted null space is defined as follows. %, see  \cite{Cohen} and \cite{Raskutti}.
\begin{definition}[Restricted Null Space (RNS)]
	For a fixed set $S \subset \{1, 2, . . . , p\}$, the design matrix $\textbf{X}$  satisfies the restricted null space property with respect to S, denoted as RN(S) if,
	\begin{align*}
	\mathbb{C}(S,1) \cap null(X) = {0}
	\end{align*}	
	where the set $\mathbb{C}(S,1)$ is defined as in Eq. (\ref{conCone}).
\end{definition}
We illustrate the RNS property using a couple of examples as follows.
\begin{example}
Consider an under-determined linear system where $p =2$ and $n =1$, but the true support is $S = \{2\}$. Let us suppose that $\textbf{X}= ( 1, 2 )$, %and we observe Y = 3 and .\\
the true coefficient vector to be estimated is of the form $\beta^0 = (0 \ b)'$, where $b \in \mathbb{R}$. The null(X) is the set of scaler multiples of $( 2 \ -1)'$. Since, for all vectors that belong to null(X) also satisfy  $|2*\alpha| \geq |- \alpha|$, hence a unique solution exists. The solution will be, where the $\ell_1$ ball intersects the null(X) translated by $(0 \ b)'$, as the error vector has to satisfy $\| \Delta_{S^c} \|_1 = |2*\alpha| \leq |- \alpha - b| = \| \Delta_S \|_1$. Therefore the solution of $\ell_1$ constrained optimization problem is  $\hat{\beta }^{\star}= (0 \ b)'$, and thus it estimates the true $\beta^0$ exactly.

\end{example}

\begin{example}
We consider another similar problem, where $p =2$, $n =1$ and the true support is $S = \{2\}$. Given $\textbf{X} = \{ 2, 1 \}$. %and Y = 3 we have to estimate true $\beta^0 = \{0, 3 \}$.\\
The true coefficient vector to be estimated is of the form $\beta^0 = (0 \ b)'$. The null(X) is the set of scaler multiples of the vector $( 1 \ -2 )'$. For all vectors that belong to null(X) does not satisfy $\| \Delta_{S^c} \|_1 = |-2 \alpha| \not> |\alpha| = \| \Delta_S \|_1$, hence the design matrix does not satisfy the RNS property. The solution of BPLP is  $\hat{\beta }^{\star}= (b/2 \ 0)'$, and thus it fails to estimate the true $\beta^0$.

\end{example}

\subsection{Sufficient Conditions for RNS}
It has been proven that the Restricted Isometry Property (RIP) and  Mutual Incoherence Property (MIP) are sufficient for the restricted nullspace property to hold.
%\subsubsection{Restricted Isometry Property(sub-matrix incoherence)}
%Restricted Isometry property(RIP) is a fundamental property which says that if a design matrix satisfies the RIP of order s(order of sparsity), then the true $\beta^0$ can be reconstructed even for high dimensional noiseless case when $p>n$.
The RIP is defined as follows, see \cite{RIP} and \cite{Hastie}.
\begin{definition}[Restricted Isometry Property (RIP)]
For a fixed $s\leq p$, the design matrix $\textbf{X}$ is said to satisfy the RIP of order $s$ with isometry constant $\delta_s$ if
\begin{align*}
(1-\delta_s) \| u\|^{2}_{2} \leq  \| X_s u\|^{2}_{2}  \leq \| u\|^{2}_{2}  (1+\delta_s) \quad \forall u \in \mathbb{R}^{s}
\end{align*} 
or equivalently 
\begin{align*}
 \frac{  \| X_s u\|^{2}_{2}}{\| u\|^{2}_{2}}  \in [(1-\delta_s),  (1+\delta_s)]
\end{align*}
In words, the matrix $X_{s}^{'}X_s$ has their eigenvalues
in $[(1-\delta_s),  (1+\delta_s)]$, for all subsets $S$ of size at most $s$.
\end{definition}
RIP is a sufficient condition for the RNS to hold, as given by the following proposition, we refer to \cite{RIP} and \cite{Hastie} for the proof.

\begin{proposition}[RIP implies RNS]
If RIP$(2s,\delta_{2s})$ holds with $\delta_{2s} < 1/3$, then the uniform RNP holds for all subsets $S$ of size at most $s$.%, therefore $\ell_1$ minimization is exact for all support vector of size at most s.
\end{proposition}
%\subsubsection{Mutual Incoherence Property}
In the following, we define Mutual Incoherence Property (MIP). %The MIP is a special case of the RIP. The MIP is a measure of the orthonormality among pairs of columns of the design matrix X, where as RIP implies constraint on much larger sub-matrices of$\textbf{X}$ to have nearly orthonormal columns. 
%RIP was previously used for exact recovery in the noiseless case, now it is shown to be sufficient for stable recovery in the noisy case also.

\begin{definition}[Mutual Incoherence Property (MIP)] 
Assuming that the columns of the design matrix $\textbf{X}$ are centred, MIP is the maximum absolute pairwise correlation, see \cite{MIP}.
\begin{align*}
	M(X) = max_{(i \neq j)} | x_{i}^{'}x_j|
\end{align*}
\end{definition}
A low pairwise incoherence is required for the exact estimation in noiseless case and for stable estimation in noisy case.

\begin{proposition}[MIP implies RNS]
Suppose that for some $s \leq p$, the MIP satisfies the bound $M(X) < \frac{1}{3s}$, then the uniform RNP of order $s$ holds. %, therefore $\ell_1$ minimization is exact for all support vector of size at most s.
\end{proposition}
For proof of the above proposition, we refer to \cite{Hastie}.
The main advantage of MIP is that it can be easily computed in $O(np^2)$ time and a major disadvantage is that it gives very conservative bound for $\ell_1$ penalized estimation. %and may fail to capture the actual performance of $\ell_1$ relaxation in practice.

\section{Constraints and Conditions }
%In this section, we discuss various constraints and conditions required for stable estimation in noisy case. We note that it is impossible to estimate exactly in noisy case. 
In this section, we study various conditions required on design matrix $\textbf{X}$ to establish oracle results for the Lasso. %for exact recovery in noiseless setting, and for stable recovery in the noisy case. 
We also describe cone constraints on error vector $\Delta$ and beta-min condition on true coefficient vector $\beta^0$.
\subsection{Convex Cone Constraints}
Here, we define a subset 
\begin{align} \label{conCone}
\mathbb{C}(S,L) = \{ \Delta \in \mathbb{R}^p| \;\| \Delta_{S^c}\|_1 \leq L \| \Delta_{S}\|_1 \}
\end{align}
for some $L \geq 1$. %This corresponds to the cone of vectors whose $\ell_1$-norm off the support is dominated by the $\ell_1$-norm on the support.
This corresponds to the cone of vectors which is a subset of $\mathbb{R}^p$, and $\Delta$ is restricted to lie in this subset.

It is easy to prove that the Lasso error $\Delta^{\star} =  \beta^{\star} - \beta^0$, for the noiseless case (Eq. \ref{noiselessLasso}) always satisfies a cone constraint with $L=1$. Suppose $ \beta^{\star}$ is a solution of the Lasso and $\beta^0$ the vector of true parameters, then due to optimality the following inequality holds.
\begin{align}
	\frac{1}{n}\| \textbf{X} \beta^{\star} - \textbf{X} \beta^0 \|_{2}^{2} +  \lambda \| {\beta}^{\star} \|_1 & \leq \frac{1}{n}\|  \textbf{X} \beta^0  - X\beta^{0} \|_{2}^{2} + \lambda \| \beta^0 \|_1 \nonumber \\
	 \frac{1}{n}\| \textbf{X} \beta^{\star} - \textbf{X} \beta^0 \|_{2}^{2} & \leq  \lambda (\| \beta^0 \|_1 - \| {\beta}^{\star} \|_1 ) \nonumber \\
	 0 & \leq  \|\beta_{S}^{\star} - \beta^0_S \|_1 -  \| {\beta}^{\star}_{S^c} \|_1 \nonumber \quad (\because \beta^0_{S^c} = 0) \\
	 & \leq   \| \Delta_S^{\star} \|_1 - \| \Delta^{\star}_{S^c} \|_1 \nonumber \\
	 \| \Delta^{\star}_{S^c} \|_1 & \leq \| \Delta^{\star}_S \|_1 \label{eq:noiselesscone}
\end{align}
Similarly, it can be shown that for the noisy case, the Lasso error $\Delta = \hat{\beta}- \beta^0$ satisfies the following cone constraint for a suitable choice of regularization parameter $\lambda$ (see section \ref{choiceL} for proof).
\begin{align}
\| \Delta_{S^c} \|_1 \leq 3 \| \Delta_{S} \|_1 
\end{align}
In general, the following theorem holds for the error vector $\Delta$.
\begin{theorem}
The Lasso error, $\Delta = \hat{\beta}- \beta^0$  is restricted to the cone constraint $\mathbb{C}(S,L)$ for some $L \geq 1$.
\end{theorem}

\subsection{Restricted Eigenvalue Condition}
The least squares objective function $ f(\beta) = \| \textbf{y} - \textbf{X} \beta \|_{2}^{2}$ is always convex in $\beta$. The function $f(\beta)$ is also strongly convex if $\nabla^2 f(\beta) = \frac{X'X}{n} \geq \gamma$ for some $\gamma > 0$, see \cite{Hastie}. For $p < n$, $|\frac{X'X}{n}| =0$, hence $f(\beta)$ can not be strongly convex. The strong convexity condition is relaxed for some subset $C \subset \mathbb{R}^p$. 
\begin{definition}[Restricted Strong Convexity]
A function $f(\beta)$ satisfies restricted $\gamma$-strong convexity at $\beta^{*}$ with respect to $C$ if 
\begin{align*}
	 \frac{  u'\nabla^2 f(\beta) u}{\| u\|^{2}_{2}}  \geq \gamma \quad \text{ for all non zero } u \in C,
\end{align*}
for all $\beta$ in neighbourhood of $\beta^*$.
\end{definition}
%In words, it is the same as lower bounding the restricted eigenvalue of the design matrix. % as given below.
%\begin{align*}
%	 \frac{ \frac{1}{n} \| \textbf{X} u\|^{2}_{2}}{\| u\|^{2}_{2}}  \geq %\gamma \quad \text{ for all non zero } u \in C
%\end{align*}
The restricted eigenvalue on the covariance matrix $\Sigma$ using the cone constraint on the Lasso error (see \cite{Bickel} and \cite{Buhlmann1}) can be given as follows.
\begin{definition}[Restricted Eigenvalue Condition (RE-condition)]
For a set $S$ with cardinality $s = |S|$ and constant $L>0$, the $(L,S,s)$-restricted eigenvalue is
\begin{align*}
	\phi^{2}(L,S,s) := min\left\lbrace { \frac{ \frac{1}{n}\| \textbf{X} \Delta \|_{2}^{2}}{\| \Delta_S \|_{2}^{2}}:   \| \Delta_{S^c}\|_1 \leq L \| \Delta_{S}\|_1 \neq 0  }  \right\rbrace.
\end{align*}
The restricted eigenvalue condition is said to be met if $\phi^{2}(L,S,s) > 0$ for all subsets $S$ of size $s$.
\end{definition}
RE-condition has been used to derive oracle results for the estimation and prediction (see \cite{Bickel}).
Since $\ell_1$ norm can be upper bounded by $\ell_2$ norm (\ref{eq:L12}), the adaptive RE condition can be defined using $\ell_2$ norm as follows.
\begin{definition}[Adaptive Restricted Eigenvalue]
For a set S with cardinality $s = |S|$ and constant $L>0$, the $(L,S,s)$-adaptive restricted eigenvalue constant is 
\begin{align*}
	\phi_{adap}^{2}(L,S,s) := min\left\lbrace { \frac{ \frac{1}{n}\| X\Delta \|_{2}^{2}}{\| \Delta_S\|_{2}^{2}}:  \| \Delta_{S^c}\|_1 \leq L \sqrt{s}\| \Delta_{S}\|_2 \neq 0  }  \right\rbrace.
\end{align*}
The adaptive restricted eigenvalue condition is said to be met if $\phi_{adap}^{2}(L,S,s) > 0$ for all subsets $S$ of size $s$.
\end{definition}
The adaptive restricted eigenvalue condition is useful in proving oracle results for the adaptive Lasso (see \cite{Zou}).

Finally we 

 a slightly stronger version of restricted eigenvalue condition as given in \cite{Raskutti}, we use it to derive lower bound for $\ell_2$ estimation error in later section.
\begin{definition}[(strong) Restricted Eigenvalue Condition]
For a set $S$ with cardinality $s = |S|$ and constant $L>0$, the $(L,S,s)$-strong restricted eigenvalue is
\begin{align*}
	\phi^{2}_{str}(L,S,s) := min\left\lbrace { \frac{ \frac{1}{n}\| \textbf{X} \Delta \|_{2}^{2}}{\| \Delta \|_{2}^{2}}:   \| \Delta_{S^c}\|_1 \leq L \| \Delta_{S}\|_1 \neq 0  }  \right\rbrace.
\end{align*}
The strong restricted eigenvalue condition is said to be met if $\phi^{2}_{str}(L,S,s) > 0$ for all subsets $S$ of size $s$.
\end{definition}
%We will discuss the relevant constraint set C in later sections.

%to be edited: The lower bound of $\nabla^2 f(x)$ implies that it is strongly convex.

%In the case of noisy observations, exact recovery of $\beta$ is no longer possible, and one goal is to obtain an estimate bβ such that the $\ell_2$-error  is well-controlled.

%We note that if the RE condition holds with parameters with $\phi(L,S,s) > 0 $, then the RNS property holds. 

%The square submatrices of size $\leq 2s$(why 2s) of the Gram matrix are necessarily positive definite. Indeed, suppose that for some $X\Delta$= 0 we have simultaneously restricted eigen value condition for set S and $X\Delta = 0$. Partition $\Delta$ in two sets: $\Delta = \Delta_S ∪ \Delta_Sc$. Without loss of generality, suppose that $\Delta_S = \Delta_Sc$ then $\phi(1,S,s) = 0$, and  a contradiction.

\subsection{Compatibility Condition}
For a fixed active set $S$ with cardinality $s = |S|$ and constant $L>0$, the $(L,S)$ restricted $\ell_1$ eigenvalue is
\begin{align} \label{eq:compCond}
	\phi_{comp}^{2}(L,S) := min \left\lbrace \frac{\frac{1}{n} \| \textbf{X} \Delta \|_{2}^{2} s}{\| \Delta_{S}\|_{1}^{2}} :  \| \Delta_{S^c}\|_1 \leq L \| \Delta_{S}\|_1 \neq 0 \right\rbrace.
\end{align}

\begin{definition}[Compatibility Condition]
The $(L,S)$ compatibility condition is said to be satisfied for the set $S$, if $\phi_{comp}(L,S) > 0$.
\end{definition}

%We note that $\| \beta_{S}\|_{1}^{2} \leq   \frac{s \frac{1}{n} \| X\beta\|_{2}^{2}} {\phi_{comp}(L,S)}$.

We note that, the RE-condition implies that the compatibility condition holds for all subset $S$ of size $s$. The compatibility condition depends on the set $S$ whereas the RE condition depends only on the cardinality $s = |S|$. It follows that the RE condition is stronger than the compatibility condition.  In later section, we derive bounds for prediction accuracy and estimation error under assumption of the least restrictive condition, the compatibility condition. % see Bickel et al. (2009)
%since $\| \beta\|_{2}^{2} \geq \| \beta_S \|_{2}^{2} \geq \| \beta_S \|_{1}^{2} /s.$.

%\begin{proposition}(RE condition implies compatibility condition)
%Suppose the RE condition holds for all subset of size s, with constant $\alpha$, then the compatibility condition holds with compatibility constant $\phi_{s}^{2} \geq \alpha$.
%\end{proposition}
%that follows $\phi_{adap}(L,S,s) \leq \phi_(L,S,s) \leq \phi_{comp}(L,S)$.

\subsection{Beta-min Condition}
Beta-min condition requires that the non-zero regression coefficients
 are sufficiently large, see \cite{Buhlmann1}. The beta-min, denoted by $|\beta_0|_{min}$ is defined as
:
\[
	|\beta_0|_{min} = min_{j \in S_0} |\beta_{j}^{0}|
\] 
In order to build an intuition, first we consider the noiseless Lasso (\ref{noiselessLasso}) problem with orthonormal design matrix, where $\frac{X'X}{n} = I_p$. Using KKT condition (see \cite{Boyd}) we get the following equality.
\begin{align*}
	\frac{2X'(X({\beta}^{\star} - \beta^0))}{n} &= -\lambda \tau^{\star} ,\quad where \quad \tau^{\star} = sign({\beta}^{\star})\\
	{\beta}^{\star} - \beta^0 &= \frac{-\lambda \tau^{\star} }{2}\\
	{\beta}^{\star} &= \beta^0 - \frac{\lambda \tau^{\star} }{2}\\
	\implies {\beta}^{\star}_j &= \left\lbrace \begin{array}{ll}
	sign(\beta^0)(|\beta^0| - \frac{\lambda}{2}) & \text{ if } |\beta^0| \geq \frac{\lambda}{2} \\
	0 &  \text{ otherwise } \end{array}	 \right. \\
	\implies {\beta}^{\star}_j &= S_{\frac{\lambda}{2}}(\beta^0)
\end{align*}
%Therefore $| \beta^0_{j}| > 0$ only when $ |(\beta_0)_j| & \geq \frac{\lambda}{2}$. 
Hence we conclude that, for the orthonormal design matrix, the Lasso (noieless case) will select a predictor if its corresponding true regression coefficient $| \beta^0_{j}| > \frac{\lambda}{2}$. For the standard Lasso the beta-min condition (see \cite{Buhlmann1}) is defined as follows.
\begin{definition}[Beta-min Condition]
The beta-min condition is said to be satisfies if the following holds
\begin{align}
	|\beta_0|_{min} \geq \frac{4 \lambda s_0}{\phi^{2}_{comp}(L,S_0)}
\end{align}
\end{definition}

In section \ref{sec:varSelection}, we will derive beta-min condition under IR-condition and we will also show that the beta-min condition is important for exact variable selection.

\subsection{Irrepresentability Condition}
The irrepresentable condition depends on the covariance of the predictors ($\Sigma$) and the signs of the unknown true parameter $\beta^0$, see \cite{Zhao}. %Lasso selects the true active set $S_0$ if and (almost) only if the predictors that are not in the $S_0$ are `irrepresentable”  by predictors that are in the true active set $S_0$, see \cite{Zhao}.
For a fixed set $S \subset \{1, 2, . . . , p\}$ and true regression coefficient $\beta$, the IR condition (weak) is given as
\[
	\| \Sigma_{21} \Sigma^{-1}_{11} sign(\beta_{S})\|_{\infty} < 1.
\]
Since we do not know $sign(\beta)$ before hand, we need the Irrepresentable Condition to hold for every possible combination of different signs and placement of zeros, so we use a modified version of IR conditions which involves only the design matrix but not the true coefficient, as given in \cite{Sara3}.

\begin{definition}[Irrepresentable Condition]
The irrepresentable condition is said to be met for the set $S$, with cardinality $s = |S|$, if the following holds:
\begin{align*}
	\|\Sigma_{12}(S) \Sigma^{-1}(S) \tau_S  \|_{\infty} < 1, \quad \forall \tau_S \in \mathbb{R}^s \; such \; that \;   
	\| \tau_S \|_{\infty} \leq 1.
\end{align*}
The weak irrepresentable condition holds for a fixed $\tau_S$ if 
\begin{align*}
	\|\Sigma_{12}(S) \Sigma^{-1}(S) \tau_S  \|_{\infty} \leq 1.
\end{align*}
For some $0<\theta <1$, the $\theta$-uniform irrepresentable condition holds if
\begin{align*}
	\mathop{max}_{\| \tau_S \|_{\infty} \leq 1} \|\Sigma_{12}(S) \Sigma^{-1}(S) \tau_S  \|_{\infty} \leq \theta
\end{align*}
\end{definition}
%It can be interpreted as, the total amount of an irrelevant covariate represented by the covariates in the true active set $S_0$ should not reach 1. When a redundant variable is strongly correlated with the predictors in the true active set $S_0$, it is impossible to distinguish it from the true predictors with any amount of data and any amount of regularization.
%To correctly identify the active set $S_0$, the irrelevant covariates should be roughly orthogonal to the relevant predictors. 
We note that the IR condition implies that false positives are not selected by the Lasso, whereas for controlling the false negatives beta-min condition is also required.

\subsection{An Example}
In this section we illustrate the various conditions defined above using a simple example. Let $S_0 = \{1,2,3,4 \}$ be the active set,  $\Sigma =  \frac{X'X}{n} $ and is given as
\begin{align*}
	\Sigma =  \left[ \begin{array}{ccccc}
	1  & 0  & 0 & 0 & \rho \\
   0   & 1  & 0 & 0 & \rho \\
   0   & 0  & 1 & 0 & \rho \\
   0   & 0  & 0 & 1	& \rho \\
 \rho & \rho & \rho& \rho & 1
	\end{array} \right],
\end{align*}
where the active variables are uncorrelated and the fifth variable is equally correlated with all active covariates. We partition  $\Sigma$ as in Eq (\ref{eq:covMat}), and we analyse for what values of  $\rho$ various conditions hold.

It is easy to check that for $\rho \geq \frac{1}{2}$, $\Sigma$ is positive semi definite.
First, we consider $\rho = \frac{1}{2}$, though the $\Sigma$ is not invertible and $\Lambda_{min}(\Sigma) = 0$, it satisfies RIP of order 4 with constant $\delta_4 = 0.866$. Since RIP implies Restricted eigenvalue(and compatibility condition), see \cite{Sara3}, therefore for $\rho = \frac{1}{2}$, the $\Sigma$ satisfies the RE and compatibility conditions. For a particular choice of $\tau_S = \{1,1,1,1 \}$, $\Sigma$
does not satisfy the IR condition.

%Now, we consider $\rho = \frac{1}{4}$, then $\Lambda_{min}(\Sigma) > 0$. The RIP of order 4 is satisfied with constant $\delta_4 = 0.433$, and therefore RE(and compatibility) condition holds but for  $\tau_S = \{1,1,1,1 \}$, $\Sigma$ does not satisfy the IR condition.
 
For $\rho < \frac{1}{s} = \frac{1}{4}$, it also satisfies the IR condition. The RIP with small constants implies the IR condition, for exact form and proof we refer to \cite{Sara3}. To check compatibility using approximations and several examples (i.e. toeplitz, block diagonal matrix etc.) where the compatibility condition holds are given in \cite{Sara3}.  %Various matrices (i.e. toeplitz, block diagonal) that satisfies the restricted eigenvalues(and compatibility condition) are discussed in \cite{Sara3}.

\section{Some Useful Bounds}
In this section, we derive important bounds which are required for proving oracle bounds for the Lasso (Eq. \ref{eq:lasso}). 
\begin{lemma}[Basic Inequality]:
\begin{align}
	\| \textbf{X} \hat{\beta} - X\beta^0 \|_{2}^{2}/n +\lambda \|\hat{\beta} \|_1 \leq  2 \epsilon^{'} X(\hat{\beta} - \beta^0)/n + \lambda \|\beta^0 \|_1 \label{eq:basicInq}
\end{align}
\end{lemma}
\begin{pf}
For an optimal solution $\hat{\beta}$ and a feasible solution $\beta^0$ the following holds.
\begin{gather*}
	\| Y - \textbf{X} \hat{\beta} \|_{2}^{2}/n +\lambda \|\hat{\beta} \|_1 \leq
	 \| Y - \textbf{X} \beta^0 \|_{2}^{2}/n +\lambda \|\beta^0 \|_1\\
	\implies \| \epsilon - \textbf{X}(\hat{\beta} - \beta^0) \|_{2}^{2}/n +\lambda \|\hat{\beta} \|_1  \leq \| \epsilon \|_{2}^{2}/n +\lambda \|\beta^0 \|_1\\
	\implies \| \textbf{X} \hat{\beta} - X\beta^0 \|_{2}^{2}/n  -  2 \epsilon^{'} X(\hat{\beta} - \beta^0)/n +\lambda \|\hat{\beta} \|_1  \leq \lambda \|\beta^0 \|_1
\end{gather*}
Rearranging the above inequality we get the basic inequality, see \cite{Buhlmann1}.
\qed
\end{pf}
The term $ 2 \epsilon^{'} X(\hat{\beta} - \beta^0)/n$ is called stochastic process part, and it can be bounded by the $\ell_1$ norm of the Lasso error. Applying Holder's inequality (\ref{eq:holder}), we get the following.
\begin{align*}
	 2 \epsilon^{'} X(\hat{\beta} - \beta^0)/n \leq
	 2 \| \epsilon^{'} X/n \|_{\infty}\; \| (\hat{\beta} - \beta^0)\|_1
\end{align*}
To overrule the stochastic process part the (regularization parameter) penalty is defined to be $\lambda \geq 2 \lambda_0$ where $\lambda_0 \geq 2 \| \epsilon^{'} X/n \|_{\infty}$ with high probability.
Therefore we can write
\begin{align}
	 2 \epsilon^{'} X(\hat{\beta} - \beta^0)/n \leq
	 \frac{\lambda}{2} \| (\hat{\beta} - \beta^0)\|_1. \label{eq:stocastic}
\end{align}

%\subsection{Bounds on \ell_2 Norms of true Parameter Vector}

\subsection{Choice of Regularization Parameter (Penalty)}
In this section, we prove that a good choice for the regularization parameter $\lambda $ is of order $\sigma \sqrt{\frac{\log p}{n}}$.\\
\\
Let $\{x_1, ..., x_p \}$ denote the columns of the design matrix $X$. Then the random variable $x^{'}_{j} \epsilon$ is distributed as $N(0, \sigma^2/n)$. We have the following inequality.
\begin{align*}
 \mathbb{P}(|x^{'}_{j} \epsilon|/n \geq t) \leq 2 e^{\frac{-nt^2}{2 \sigma^2}}
\end{align*}

As  $\| \epsilon^{'} X/n \|_{\infty}$ corresponds to the maximum over $p$ such variables, the union bound gives the following inequality.
\begin{align*}
	\mathbb{P}( \| \epsilon^{'} X/n \|_{\infty} \geq t) & \leq 2 p e^{\frac{-nt^2}{2 \sigma^2}}\\
	& = 2 e^{\frac{-nt^2}{2 \sigma^2}+ \log p}
\end{align*}
By setting $ t = \sigma \sqrt{\frac{\tau \log p}{n}}$ for some $\tau > 2$,  we get $\lambda = const. \sigma \sqrt{\frac{\log p}{n}}$ with high probability. 
For rest of the article we choose $\lambda$ to be the order of $\sigma \sqrt{\frac{\log p}{n}}$.

\subsection{Choice of the Parameter L} \label{choiceL}
In this section, we prove that a good choice for the parameter $L$, which is used for cone constraint, is of order $L = \frac{c+1}{c-1}$, where $c$ is the defined as follows.
\begin{align*}
	\lambda & > 2 \| \epsilon^{'} X/n \|_{\infty}\\
 \lambda & = c (2 \| \epsilon^{'} X/n \|_{\infty})\\
 \frac{\lambda}{c} &= 2 \| \epsilon^{'} X/n \|_{\infty} 
 \end{align*}
%In general , c is chosen as $c=2$ , therefore L=3.\\
\begin{pf} %For suitable choice of $\lambda$, we can prove L = 3.\\
Using the basic inequality, and let $\Delta = \hat{\beta} - \beta^0$
\begin{align*}
	\| \textbf{X} \Delta \|_{2}^{2}/n +\lambda \|\hat{\beta} \|_1 &\leq   \frac{\lambda}{c} \|\Delta \|_1 + \lambda \|\beta^0 \|_1 \\
	0 & \leq  \frac{\lambda}{c} \|\Delta \|_1 + \lambda \| \hat{\beta }_S- \beta^0 \|_1 - \lambda \| \hat{\beta }_{S^c} \|_1 \\
	0 & \leq  \frac{\lambda}{c} \|\Delta_S  \|_1 + \frac{\lambda}{c} \| \hat{\Delta^c} \|_1 + \lambda \| \Delta_S \|_1  - \lambda \| \Delta_{S^c} \|_1 \\
	0 & \leq  \frac{\lambda (c+1)}{c} \|\hat{\beta }_S- \beta^0  \|_1  - \frac{\lambda(c-1)}{c} \| \hat{\beta }_{S^c} \|_1\\
	\| \Delta_{S^c} \|_1 &\leq \frac{c+1}{c-1} \| \Delta_S  \|_1\\
	\text{Or, equivalently}&\\
	\|\beta_{S^c}  \|_1 & \leq L \| \hat{\beta}_S - \beta^0  \|_1
\end{align*}
We proved that $  L =  \frac{c+1}{c-1}$. For rest of the article we choose $c=2$, therefore $L=3$.
\qed
\end{pf}

\section{Oracle Bounds for the Lasso}
Usually, we assess quality of the Lasso estimates by measuring prediction and estimation accuracy, and variable selection consistency. In this section, we discuss various bounds for them, see  \cite{Sara1},\cite{Bickel}, \cite{Bunea}, \cite{Nicolai}, \cite{Sara2}, \cite{Sara3}),  and for a book length discussion we refer to \cite{Buhlmann1}. 
\subsection{Bounds on Prediction Loss}
Prediction loss is also known as mean squared prediction error, and it is defined as (see \cite{Hastie}):
\begin{align}
	L_{pred}(\hat{\beta};\beta^0) = \frac{1}{n}\| \textbf{X} \hat{\beta} - \textbf{X} \beta^{0} \|_{2}^{2} =  \frac{1}{n}\| \textbf{X} \Delta \|_{2}^{2}
\end{align}
We start with deriving prediction error bound for the noiseless Lasso for building understanding, then it becomes easy to derive the bounds for the noisy case. 

\subsubsection{Bounds on Prediction Loss in the noiseless case}
First, we consider the noiseless Lasso problem and we derive oracle result for the prediction without assuming any condition on the design matrix. 
\begin{theorem}
If ${\beta}^{\star}$ is an optimal solution for the noiseless Lasso problem as given in Eq. (\ref{noiselessLasso}), then it satisfies the following bound.
\begin{align} \label{eq:slow1}
	 \frac{1}{n}\| \textbf{X} ({\beta}^{\star} - \beta^0) \|_{2}^{2} & \leq  \lambda \| \beta^0 \|_1 %\approx \| \beta^0 \|_1 \sqrt{\frac{\log p}{n}}
\end{align}
\end{theorem}

\begin{pf}
If ${\beta}^{\star}$ is an optimal solution for the Eq. (\ref{noiselessLasso}), then the following basic inequality holds.
\begin{align}
	%\frac{1}{n}\| Y - \textbf{X} \hat{\beta} \|_{2}^{2} +  \lambda \| \hat{\beta}\|_1 & \leq \frac{1}{n}\| Y - X\beta^{0} \|_{2}^{2} + \lambda \| \beta^0 \|_1\\
	\frac{1}{n}\|  \textbf{X} {\Delta}^{\star} \|_{2}^{2} + \lambda \| {\beta}^{\star} \|_1 & \leq  \lambda \| \beta^0  \|_1 \\
	 \frac{1}{n}\| X {\beta}^{\star} - X\beta^{0} \|_{2}^{2} & \leq  \lambda \| \beta^0 \|_1 - \lambda \| {\beta}^{\star} \|_1 \\
	 0 & \leq  \lambda \| \beta^0 \|_1 - \lambda \| {\beta}^{\star} \|_1 \nonumber \\
	  \| {\beta}^{\star} \|_1 & \leq \| \beta^0 \|_1  \\
	  \therefore  \frac{1}{n}\| \textbf{X} {\Delta}^{\star} \|_{2}^{2} & \leq \lambda \| \beta^0 \|_1  
\end{align}
\end{pf}
Now we exploit compatibility constant to derive fast rate bounds for the prediction error. 
\begin{theorem}
If ${\beta}^{\star}$ is an optimal solution for the Eq. (\ref{noiselessLasso}) and $\phi(S) > 0$ over a fixed set S, then ${\beta}^{\star}$  satisfies the following bounds.
\begin{align}
\frac{1}{n}\|\textbf{X} {\Delta}^{\star} \|_{2} & \leq \frac{\lambda^2 s}{ \phi^{2}_{comp}}\\
\| \Delta^{\star} \|_1 &\leq \frac{2 \lambda s}{ \phi^{2}_{comp}}
\end{align}
\end{theorem}

\begin{pf}
Once again we start from the basic inequality.
\begin{align}
	 \frac{1}{n}\|  \textbf{X} {\Delta}^{\star} \|_{2}^{2} & \leq  \lambda \| \beta^0 \|_1 - \lambda \| {\beta}^{\star} \|_1 \\
	 & \leq \lambda \|{\beta_S}^{\star} - \beta^0 \|_1 - \lambda \| {\beta}^{\star}_{S^c} \|_1\\
	 \frac{1}{n}\| \textbf{X} {\Delta}^{\star} \|_{2}^{2}  & \leq \lambda \| \Delta^{\star}_S \|_1 - \lambda \| \Delta_{S^c}^{\star} \|_1 \label{eq:losslessBasic1}
	% 0 & \leq \lambda \| \Delta^{\star}_S \|_1 - \lambda \| \Delta^{\star}_{S^c} \|_1\\
	 %\implies \| \Delta^{\star}_{S^c} \|_1 & \leq \| \Delta^{\star}_S \|_1 
\end{align}
From Eq. (\ref{eq:losslessBasic1}), we can also have 
\begin{align*}
 \frac{1}{n}\| \textbf{X} {\Delta}^{\star} \|_{2}^{2} +  \| \Delta^{\star}_{S^c} \|_1 &\leq \| \Delta^{\star}_S \|_1 \\
%\because  \| \Delta_{S^c} \|_1  & \leq \| \Delta^{\star}_S \|_1 \\
\therefore \frac{1}{n}\| \textbf{X} {\Delta}^{\star} \|_{2}^{2}  & \leq \| \Delta^{\star}_{S} \|_1
\end{align*}
Substituting value for   $\| \Delta^{\star}_S \|_1$ from Eq. (\ref{eq:compCond})
\begin{align*}
	  \frac{1}{n}\| \textbf{X} {\Delta}^{\star} \|_{2}^{2} & \leq \frac{\lambda \|\textbf{X} {\Delta}^{\star} \| \sqrt{s} }{\sqrt{n} \phi_{comp}(S,L)}\\
	  \frac{1}{n}\|\textbf{X} {\Delta}^{\star} \|_{2} & \leq \frac{\lambda^2 s}{ \phi^{2}_{comp}} %\leq   \approx \frac{ s \log p}{ n \phi^{2}_{comp}}
\end{align*}
%The above bound converges faster than the bound given in (\ref{eq:slow1}).\\
Continuing from Eq. (\ref{eq:losslessBasic1}), adding $ \lambda\| \Delta^{\star} \|_1 $ to both the sides we get  
\begin{align*}
	 \frac{1}{n}\| \textbf{X} \Delta^{\star} \|_{2}^{2} + \lambda \|\Delta^{\star} \|_1 & \leq 
	 \lambda \| \Delta^{\star}_S \|_1 - \lambda \| \Delta^{\star}_{S^c} \|_1 + \lambda( \| \Delta^{\star}_S \|_1 + |\| \Delta^{\star}_{S^c}\|_1)\\
	 & \leq 2\lambda \| \Delta^{\star}_S \|_1 .
\end{align*}
Substituting value for   $\| \Delta^{\star}_S \|_1$ from Eq. (\ref{eq:compCond}) and multiplying by 2, we get
\begin{align*}
	  \frac{2}{n}\| \textbf{X} \Delta^{\star} \|_{2}^{2} + 2\lambda \|\Delta \|_1 & \leq \frac{4\lambda \|  \textbf{X} \Delta^{\star} \| \sqrt{s} }{\sqrt{n} \phi_{comp}(S,L)}.
\end{align*}
Using the inequality on r.h.s    $4ab \leq a^2 +4b^2 $.
\begin{align*}
	 \frac{2}{n}\| \textbf{X} \Delta^{\star} \|_{2}^{2} + 2\lambda \|\Delta \|_1  & \leq \frac{ \|  \textbf{X} \Delta^{\star}\|_2^{2} }{  n }+ \frac{4 \lambda^2 s}{ \phi^{2}_{comp}}\\
	  \frac{1}{n}\| \textbf{X} \Delta^{\star} \|_{2}^{2} + 2 \lambda \|\Delta^{\star} \|_1  & \leq \frac{4\lambda^2 s}{ \phi^{2}_{comp}} \\
	  \implies  2 \lambda \|\Delta^{\star} \|_1  & \leq \frac{4\lambda^2 s}{ \phi^{2}_{comp}}
\end{align*}
Therefore  $\ell_1$- estimation error bound $\| \hat{\beta} - \beta^0 \|_1 $ is
\[
  \| \Delta^{\star} \|_1 \leq \frac{2 \lambda s}{ \phi^{2}_{comp}}
 \]
 \qed
\end{pf}

\subsubsection{Bounds on Prediction Loss in the Noisy Case}
It is known that for $p<n$ and the design matrix $\textbf{X}$ with full column rank, the prediction error for the Ordinary Least Square (OLS) estimate follows $\chi^{2}_{p}$ (see \cite{Seber}):
\begin{align*}
	\| X\hat{\beta}_{OLS} - X\beta^{0} \|_{2}^{2} \sigma^2 & \sim \chi^{2}_{p} 
\end{align*}
It follows that,
\begin{align*}
	 \mathbb{E} {\frac{1}{n}\| \textbf{X} \hat{\beta}_{OLS} - X\beta^{0} \|_{2}^{2}} &= \frac{\sigma^2}{n} p
\end{align*}
When p is fixed and $n \rightarrow \infty$, the above bound tends to $0$ in probability.

We can not have similar bound for the Lasso, since for the situation where $p>n$ and $p \geq n \rightarrow \infty $, the term $\frac{p}{n}$ does not converge. Under sparsity assumption, we derive the slow rate (A) and fast rate (B) bounds for the Lasso estimator as follows. For the noisy situation, a lower bound on the regularization parameter $\lambda$ is required, that is  $\lambda \geq 2\lambda_0$, to overrule the stochastic process part.

\begin{theorem}(Prediction error bounds)\\
(A) An optimal solution $\hat{\beta} $ of the Lasso problem satisfies the following bound
\begin{align*}
	\frac{1}{n}\| X\hat{\beta} - X\beta^{0} \|_{2}^{2} &\leq
	\frac{3 \lambda}{2}\|\beta^0 \|_1 = 	const. \sigma \| \beta^0 \|_1 \sqrt{\frac{\log p  }{n} }
\end{align*}
(B) If the design matrix satisfies $\phi_{comp}(3,S)$ compatibility condition over a fix set $S$ , then an optimal solution $\hat{\beta} $ satisfies the following bound
\begin{align*}
	\frac{1}{n}\| X\hat{\beta} - X\beta^{0} \|_{2}^{2} &\leq
	\frac{9 \lambda^2 s}{4\phi^2}   = 	const. \frac{ \sigma^2 \log p}{n}  \frac{s}{\phi^2_{comp}}
\end{align*}
\end{theorem}
%\frac{1}{n}\| X\hat{\beta} - X\beta^{0} \|_{2}^{2} &\leq 	const. \frac{ \sigma^2 \log p}{n} s_0\\
%\frac{1}{n}\| X\hat{\beta} - X\beta^{0} \|_{2}^{2} &\leq 	\frac{4 \lambda^2 s_0 }{\phi_0}
Without assuming any restrictive condition on design matrix, we can derive bounds for the prediction error (A). (A) is known as slow rate since prediction error is inversely proportional to the square root of the number of observations, $n$, whereas the second condition (B) is known as fast rate because it is inversely proportion to $n$. As $n$ grows the rate of convergence for (B) will be faster. 

\begin{pf}
(A): %$\frac{1}{n}\| X\hat{\beta} - X\beta^{0} \|_{2}^{2} \leq 	6 \lambda \|\beta^0 \|_1$\\
From basic inequality (\ref{eq:basicInq}) we can get the following:
\begin{align*}
	\frac{1}{n}\| X\hat{\beta} - X\beta^{0} \|_{2}^{2}  & \leq  \frac{\lambda}{2} \|\Delta \|_1 + \lambda (\| \beta^0 \|_1- \| \hat{\beta } \|_1) \\
	0 & \leq  \frac{\lambda}{2} \|\Delta \|_1 + \lambda (\| \beta^0 \|_1- \| \hat{\beta } \|_1).
\end{align*}
Applying the inequality $  \| \hat{\beta } - \beta^0 \|_1 \leq \|\hat{\beta} \|_1 + \| \beta^0 \|_1 $, we get
\begin{gather*}
	0  \leq  \frac{\lambda}{2} (\|\hat{\beta} \|_1 + \| \beta^0 \|_1) + \lambda (\|\beta^0 \|_1- \|\hat{\beta } \|_1) \\
	\implies  \| \hat{\beta } \|_1 \leq 3 \| \beta^0 \|_1.
\end{gather*}
Considering once again the basic inequality (\ref{eq:basicInq})
\begin{align*}
	\frac{1}{n}\| X\hat{\beta} - X\beta^{0} \|_{2}^{2}  & \leq  \frac{\lambda}{2} (\|\hat{\beta} \|_1 + \| \beta^0 \|_1) + \lambda (\|\beta^0 \|_1- \|\hat{\beta } \|_1)  \\
	& \leq \frac{3 \lambda}{2} \| \beta^0 \|_1 - \frac{\lambda}{2} \| \hat{\beta } \|_1\\
	\frac{1}{n}\| X\hat{\beta} - X\beta^{0} \|_{2}^{2}  +  \frac{\lambda}{2} \| \hat{\beta } \|_1 & \leq \frac{3}{2} \| \beta^0 \|_1 \\
	\therefore \frac{1}{n}\| X\hat{\beta} - X\beta^{0} \|_{2}^{2} &\leq
	\frac{3 \lambda}{2}\|\beta^0 \|_1\\
	 & = 	const. \sigma \| \beta^0 \|_1 \sqrt{\frac{\log p  }{n} } \quad \left(\because \lambda = O(\sigma \sqrt{\frac{\log p  }{n}} \right)
\end{align*}
We note that, consistency for the prediction can be achieved only if $\| \beta^0 \|_1 \ll \sqrt{\frac{n}{\log p}}$.
\qed
\end{pf}

\begin{pf}(B)  We start with the basic inequality (\ref{eq:basicInq})
\begin{align*}
	\frac{1}{n}\| \textbf{X} \Delta \|_{2}^{2} & \leq  \frac{\lambda}{2} \| \Delta \|_1 + \lambda \|\beta^0 \|_1 - \lambda \|\hat{\beta} \|_1 \\
	& \leq \frac{\lambda}{2} \| \Delta_S+\Delta_{S^c}  \|_1 + \lambda \| \Delta_S \|_1 - \lambda \| \Delta_{S^c} \|_1
\end{align*}
 multiplying by 2 to both the sides we get
\begin{align}
	\frac{2}{n}\| \textbf{X} \Delta \|_{2}^{2} 	& \leq  3 \lambda \| \Delta_S \|_1 - \lambda \| \Delta_{S^c} \|_1. \label{eq:lassoPredErr}
\end{align}
Simplifying further
\begin{align*}
 \frac{2}{n}\| \textbf{X} \Delta \|_{2}^{2} + \lambda \| \Delta_{S^c} \|_1  &\leq 3 \lambda \| \Delta_S \|_1\\
 \therefore \frac{1}{n}\| \textbf{X} \Delta \|_{2}^{2} &\leq \frac{3}{2} \lambda \| \Delta_S \|_1\\
\end{align*}
Substituting $ \| \Delta_S \|_1 \leq  \frac{\| X\Delta \|_2 \sqrt{s}}{\sqrt{n} \phi_{comp}(3,S)}$ we get
\begin{align*}
 \frac{1}{n}\| \textbf{X} \Delta \|_{2}^{2} &\leq \frac{3 \lambda}{2}   \frac{\| X\Delta \|_2 \sqrt{s}}{\sqrt{n} \phi_{comp}(3,S)}\\
 \frac{1}{\sqrt{n}} \| \textbf{X} \Delta \|_{2} &\leq \frac{3 \lambda \sqrt{s}}{2  \phi_{comp}(3,S)} \\
 \therefore \frac{1}{n} \| \textbf{X} \Delta \|_{2}^{2} &\leq \frac{9 \lambda^2 s}{4  \phi^{2}_{comp}(3,S)}. 
\end{align*}
\qed
\end{pf}
Continuing from Eq.(\ref{eq:lassoPredErr}) and adding $ \lambda\| \Delta \|_1 $ both the sides we get
\begin{align*}
	\frac{2}{n}\| \textbf{X} \Delta \|_{2}^{2} + \lambda\| \Delta \|_1 &\leq
	 3 \lambda \| \Delta_S \|_1 - \lambda \| \Delta_{S^c} \|_1 + \lambda(\| \Delta_S  \|_1 + \| \Delta_{S^c} \|_1) \\
	& \leq 4 \lambda \| \Delta_S \|_1.
\end{align*}
Substituting $ \| \Delta_S \|_1 \leq  \frac{\| X\Delta \|_2 \sqrt{s}}{\sqrt{n} \phi_{comp}(3,S)}$
\begin{align*}
	\frac{2}{n}\| \textbf{X} \Delta \|_{2}^{2} + \lambda\| \Delta \|_1 &\leq
	4 \lambda  \frac{\| X\Delta \|_2 \sqrt{s}}{\sqrt{n}\phi_{comp}}
\end{align*}
Applying the inequality $4ab \leq a^2 +4b^2$\\
\begin{align*}
	\frac{2}{n}\| \textbf{X} \Delta \|^{2}_{2} + \lambda\| \Delta \|_1 & \leq  \frac{1}{n}\| \textbf{X} \Delta \|^{2}_{2} + \frac{4 s \lambda^2}{\phi^{2}_{comp}}  \\
	\frac{1}{n}\| \textbf{X} \Delta \|^{2}_{2} + \lambda\| \Delta \|_1 & \leq  \frac{4 s \lambda^2}{\phi^{2}_{comp}}
\end{align*}
Basically the above inequality gives two bounds.\\
(1) The prediction error bound \\
\begin{align*}
	\frac{1}{n}\| \textbf{X} \Delta \|^{2}_{2}  & \leq  \frac{4 s \lambda^2}{\phi^{2}_{comp}}
\end{align*}
(2) $\ell_1$ - estimation error bound.
\begin{align}
	\lambda\| \Delta \|_1 & \leq  \frac{4 s \lambda^2}{\phi^{2}_{comp}} \nonumber \\
	\| \hat{\beta} - \beta^0 \|_1 & \leq  \frac{4 s \lambda}{\phi^{2}_{comp}} \label{noiseL1Err}
\end{align}
The above result (and proof) is from \cite{Buhlmann1}.%, and the derivation is essentially the same as theirs.

\subsection{Bounds on Parameter Estimation Loss}
Parameter estimation loss is also known as $\ell_2-$error and is defined as the $\ell_2$-norm loss between a Lasso estimated $\hat{\beta}$ and the true regression vector $\beta^{0}$.
\begin{align}
	L_{2}(\hat{\beta};\beta^0) = \| \hat{\beta} - \beta^{0} \|_{2}^{2}
\end{align}
\\
We have already computed bound for the $\ell_1$ estimation error Eq. (\ref{noiseL1Err}).
\begin{align*}
	\| \hat{\beta} - \beta^0 \|_1 & \leq  \frac{4 s \lambda}{\phi^{2}_{comp}}
\end{align*}
 Different bounds on estimation error with similar scaling have been discussed in various papers, i.e., \cite{Bickel}, \cite{Bunea}, \cite{Candes} and \cite{Nicolai}.

\subsubsection{$\ell_2$ Estimation Bound for Noiseless Lasso}
Now we obtain optimal bound on $\ell_2$-error under assumption of more restrictive condition, so called the (strong) restricted eigenvalue condition.
\begin{theorem}
For a fixed set S, and with $L> 0$ and $\phi_{str}(L,S,s) > 0$, the solution of the noiseless Lasso Eq.(\ref{noiselessLasso}), ${\beta}^{\star}$, satisfies the following bound.
\begin{align*}
	 \| {\beta}^{\star} - \beta^{0} \|^{2}_{2} \leq \frac{ s \lambda^2}{\phi^4_{str}} 
\end{align*}
\end{theorem}

\begin{pf}
We start with the simplified basic inequality for the noiseless Lasso as in Eq. (\ref{eq:losslessBasic1}) and atfer adding $\lambda \| \Delta^{\star}_{S^c} \|_1$ to both the sides, we get
\begin{align*}
 \frac{1}{n}\| \textbf{X} \Delta^{\star} \|_{2}^{2} + \lambda \|\Delta^{\star}_{S^c} \|_1 & \leq  \lambda \| \Delta^{\star} \|_1 \\
 \therefore \frac{1}{n}\| \textbf{X} \Delta^{\star} \|_{2}^{2}  & \leq  \lambda \| \Delta^{\star} \|_1 .
\end{align*}
Applying condition from Eq. (\ref{eq:L12}) we get
\[
	\frac{1}{n}\| \textbf{X} \Delta^{\star}  \|_{2}^{2}  \leq \sqrt{s} \lambda \| \Delta^{\star} \|_2.
\]
After applying (strong) restricted eigenvalue condition
$\frac{\| \textbf{X} \Delta^{\star} \|_{2}^{2} }{n} \geq \phi^{2}_{str} \| \Delta^{\star} \|_{2}^{2}$,
we get the following.
\begin{align*}
	\phi^2_{str} \| \Delta^{\star} \|_{2}^{2} & \leq    \sqrt{s} \lambda  \| \Delta^{\star} \|_2\\
	\therefore \| \Delta^{\star} \|_{2}^{2} & \leq  \frac{ s \lambda^2}{ \phi^4_{str}} 
\end{align*}
\qed
\end{pf}
 
 \subsubsection{$\ell_2$ Estimation Bound for Noisy Case}
In the following, we show the optimal bound on $\ell_2$-error for the noisy case.
\begin{theorem}
For a fixed set $S$, and with $L > 1$, $\phi^{2}_{str}(L,S,s) > 0$ and $\lambda \geq 2 \lambda_0 $, if $\hat{\beta}$ is a solution of the Lasso problem defined in (\ref{eq:lasso}), then $\hat{\beta}$  satisfies the following inequality.
\begin{align*}
	 \| \hat{\beta} - \beta^{0} \|^{2}_{2} \leq \frac{9 \lambda^2 s}{4 \phi^4_{str}} 
\end{align*}
\end{theorem}

\begin{pf}
We start with the simplified basic inequality as in Eq.(\ref{eq:lassoPredErr}) and after adding $\lambda \| \Delta_{S^c} \|_1$
to both the sides, we get
\begin{align*}
	\frac{1}{n}\| \textbf{X} \Delta \|_{2}^{2}  + 2\lambda \| \Delta_{S^c} \|_1 & \leq \frac{ \lambda }{2} \| \Delta\|_1 + \lambda \| \Delta_S \|_1+\lambda \| \Delta_{S^c} \|_1\\
	\frac{1}{n}\|\textbf{X} \Delta \|_{2}^{2}  + 2\lambda \| \Delta_{S^c} \|_1 & \leq \frac{ 3\lambda }{2} \| \Delta\|_1	\\
	\therefore \frac{1}{n}\|\textbf{X} \Delta \|_{2}^{2}   & \leq \frac{ 3\lambda }{2} \| \Delta\|_1	.
\end{align*}
Using inequality from Eq.(\ref{eq:L12}) we get
\[
	\frac{1}{n}\| \textbf{X} \Delta  \|_{2}^{2}  \leq \frac{3 \lambda }{2} \sqrt{s} \| \Delta\|_2.
\]
After applying (strong) restricted eigenvalue condition 
$\frac{\| \textbf{X} \Delta \|_{2}^{2} }{n} \geq \phi^{2}_{str} \| \Delta \|_{2}^{2}$,
we get 
\begin{align*}
	\phi^2_{str} \| \Delta \|_{2}^{2} & \leq   \frac{3 \lambda }{2} \sqrt{s} \| \Delta\|_2\\
	\therefore \| \Delta \|_{2}^{2} & \leq  \frac{9 \lambda^2 s}{4 \phi^4_{str}} 
\end{align*}
\qed
\end{pf}

The oracle results for prediction and estimation are summarized in the Table \ref{table:lasso}.
\begin{center}
\begin{table}[h!] 
		%\begin{tabular}{| p{5mm}p{10mm}p{25mm}p{25mm}p{25mm}p{25mm}|}
		\begin{tabu} {| l  l  |c | c|}
		\hline 
		 Property &  Design Condition & \multicolumn{2}{|c|}{Lasso Method} \\
		 \hline 
		 &&Noiseless&Gaussian Noise\\  \hline
		Pred. Error (slow rate) & No requirement & $\lambda \| \beta^0 \|_1$ & $\frac{3 \lambda}{2}\|\beta^0 \|_1 $ \\	\hline 	
		 Pred. Error (fast rate) & Compatibility & $\frac{\lambda^2 s}{ \phi^{2}_{comp}} $&$
	\frac{9 \lambda^2 s}{4\phi^2_{comp}}$ \\		\hline 
		Est. Error ($\ell_1$ norm) & Compatibility & $\frac{2 \lambda s}{ \phi^{2}_{comp}} $ & $\frac{4 \lambda s}{\phi^{2}_{comp}} $\\	\hline 
		 Est. Error ($\ell_2^2$ norm) & (strong) RE & $\frac{  \lambda^2 s}{ \phi^4_{str}}  $ & $\frac{9 \lambda^2 s}{4 \phi^4_{str}} $ \\	\hline
		 %Variable Selection & Irrepresentable condition& ${S}^{\star} \subset S $ & $\hat{S} \subset S $ \\	\hline
		\end{tabu}	
		\caption{Oracle results and sufficient conditions required for prediction and estimation}\label{table:lasso}
\end{table}
\end{center}

\subsection{Variable Selection} \label{sec:varSelection}
Variable selection is also known as support recovery, and the following loss function is used:
\begin{align}
	L_{vs}(\hat{\beta};\beta^0) = \left\lbrace \begin{array}{ll}
	0 &if sign(\hat{\beta}_i) = sign(\beta^{0}_{i}) \\
	1 & otherwise
\end{array}	 \right.
\end{align}
This checks whether the Lasso estimated $\hat{\beta}$ and the true regression vector $\beta^{0}$ have the same sign support. 
It has been proven that (see \citep{Buhlmann1}) the Lasso does not select the true active set $S_0$, unless we make the following assumptions:
\begin{enumerate}
\item Irrepresentable Condition
\[ \| \Sigma_{21} \Sigma^{-1}_{11} sign(\beta^0)\|_{\infty} \leq 1 \]
\item $\lambda \geq  2 \lambda_0$
\item beta-min condition $ |\beta^0|_{min} \gg \lambda $
\end{enumerate}
%\subsubsection{Variable Selection Consistency for Noiseless case} 
The last two conditions are required to address the false negatives.  %Keeping aside the beta-min condition, we mainly focus on avoiding false positives. 
In the next two theorems, we show the variable selection consistency for noiseless case, that is the irrepresentable condition is  (almost) necessary and sufficient condition for the Lasso to select only variables in the true active set. These results and proofs are from \cite{Sara3}.

\begin{theorem}[IR is sufficient condition for variable selection]
Let ${S}^{\star}$ be the variables selected by the Lasso (\ref{noiselessLasso}). Suppose the irrepresentable condition is met for a fixed active set $S$, then ${S}^{\star} \subset S$ and 
\begin{align*}
	\| {\beta}^{\star}_{S} - \beta^0 \|_{\infty} \leq \lambda  \mathop{sup}_{\| \tau_{S}\|_{\infty \leq 1}} \| \Sigma^{-1}_{11} \tau_S \|_{\infty}/2	
\end{align*}
\end{theorem}
\begin{pf}
We consider the noiseless case Eq.(\ref{noiselessLasso}). 
Compute the KKT condition as follows.
\begin{align*}
	\frac{1}{n}2X' (X\beta - X\beta^0)+\lambda \tau = 0 ,
\end{align*} 
where $\| \tau \|_{\infty} \leq 1$ and $ \tau'\beta= \| \beta \|_1 $. 
Let ${\beta}^{\star}, {\tau}^{\star}$ be a solution of the Lasso, then we get,
\begin{align*}
	\frac{X'X}{n}2 ( {\beta }^{\star}- \beta^0) &= - \lambda {\tau^{\star} }\\
	\implies 2\Sigma ({\beta }^{\star}- \beta^0)& = - \lambda {\tau}^{\star}.
\end{align*}
Without loss of generality we can assume that the first $s$ variables are the active variables, and we partition the $\Sigma$ as in (\ref{eq:covMat}), ${\beta }^{\star}= (\beta_1 \ \beta_2)' $ and ${ \tau }^{\star}= (\tau_1 \ \tau_2)' $ accordingly (we omit the star and superscript for simplicity). We get the following equation.
\begin{align*}
\Sigma = \left[ \begin{array}{cc}
	\Sigma_{11} & \Sigma_{12}\\
	\Sigma_{21} &\Sigma_{22}	
	\end{array} \right] 	\left( \begin{array}{c} \beta_1- \beta^{0}_{1} \\ \beta_2- \beta^{0}_{2} \end{array} \right) & = \left( \begin{array}{c} \tau_1 \\ \tau_2 \end{array} \right)
\end{align*}
We note that $\| \beta^{0}_{2} \|_1 = \| \beta^{0}_{S^c}\|_1 = 0 $.
We get the following two equations after simplification:
\begin{align*}
	2\Sigma_{11}(\beta_1-\beta^{0}_{1}) + 2 \Sigma_{12} \beta_2 &= -\lambda \tau_1\\
	2\Sigma_{21}(\beta_1-\beta^{0}_{1}) + 2 \Sigma_{22} \beta_2 &= -\lambda \tau_2
\end{align*}
Multiplying the first equation with $\Sigma^{-1}_{11}$ we get
\begin{align*}
	2(\beta_1-\beta^{0}_{1}) + 2\Sigma^{-1}_{11} \Sigma_{12} \beta_2 &= -\lambda \Sigma^{-1}_{11} \tau_1\\
	2\Sigma_{21}(\beta_1- \beta^{0}_{1}) + 2 \Sigma_{22} \beta_2 &= -\lambda \tau_2.
\end{align*}
Multiplying the first equality by$-\beta^{'}_{2} \Sigma_{21}$ and the second by $-\beta^{'}_{2}$, we obtain the following.
\begin{align*}
	-2\beta^{'}_{2} \Sigma_{21}(\beta_1-\beta^{0}_{1}) - 2\beta^{'}_{2} \Sigma_{21} \Sigma^{-1}_{11} \Sigma_{12} \beta_2 &= +\lambda \beta^{'}_{2} \Sigma_{21} \Sigma^{-1}_{11} \tau_1\\
	-2\beta^{'}_{2} \Sigma_{21}(\beta_1- \beta^{0}_{1}) - 2 \beta^{'}_{2}  \Sigma_{22} \beta_2 &= +\lambda \beta^{'}_{2}\tau_2 = \lambda\| {\beta}^{\star}_{S^c}\|_1
\end{align*}
Subtracting the first equality with the second one we get
\begin{align*}
	 2\beta^{'}_{2}(  \Sigma_{22} - \Sigma_{21} \Sigma^{-1}_{11} \Sigma_{12} )\beta_2 &= \lambda \beta^{'}_{2} \Sigma_{21} \Sigma^{-1}_{11} \tau_1 - \lambda\| {\beta}^{\star}_{S^c}\|_1.
\end{align*}
Using Holder inequality for the following expression.
\begin{align*}
\beta^{'}_{2} \Sigma_{21} \Sigma^{-1}_{11} \tau_1 &\leq \lambda \|\beta_2 \|_1 \| \Sigma_{21} \Sigma^{-1}_{11} \tau_1 \|_{\infty}\\
&\leq \|{\beta}^{\star}_{S^c} \|_1
\end{align*} 
It follows that
\begin{align*}
	 2\beta^{'}_{2}(  \Sigma_{22} - \Sigma_{21} \Sigma^{-1}_{11} \Sigma_{12} )\beta_2 \leq 0.
\end{align*}
Since $(  \Sigma_{22} - \Sigma_{21} \Sigma^{-1}_{11} \Sigma_{12} )$ is a positive definite matrix it is a contradiction. Hence $\beta_{2}  = 0 \implies {\beta^{\star}_{S^c}}  = 0$. Therefore ${S}^{\star} \subset S_0$, irrepresentable condition implies no false positive selection.
\\
\\
Under irrepresentable condition, $\beta_{2}  = 0 \ ( {\beta^{\star}_{S^c}}  = 0)$, we get the following equation using the KKT condition.
\begin{align*}
	\Sigma = \left[ \begin{array}{cc}
	\Sigma_{11} & \Sigma_{12}\\
	\Sigma_{21} &\Sigma_{22})	
	\end{array} \right] 	 \left( \begin{array}{c} \beta_1- \beta^{0}_{1} \\ 0 \end{array}\right) & = \left( \begin{array}{c} \tau_1 \\ \tau_2 \end{array}\right) 
\end{align*}
After simplification we get
\begin{align*}
	2\Sigma_{11}(\beta_1-\beta^{0}_{1})  &= -\lambda \tau_1\\
	2\Sigma_{21}(\beta_1-\beta^{0}_{1}) &= -\lambda \tau_2.
\end{align*}
From first equality, it follows that
\begin{align*}
	(\beta_1-\beta^{0}_{1}) &= - \lambda \Sigma^{-1}_{11} \tau_1/2\\
	 \implies  \| (\beta_S-\beta^{0}) \|_{\infty} &\leq \lambda \| \Sigma^{-1}_{11} \tau_S \|_\infty /2.
\end{align*}
Suppose for some $j \in S$, ${\beta}^{\star}_j = 0 $ and $\beta^{0}_{j} > \lambda \| \Sigma^{-1}_{11} \tau_S \|_\infty /2 $, then we would have
\[ 
	 \| \beta_S -\beta^{0} \|_{\infty} >  \lambda \| \Sigma^{-1}_{11} \tau_S \|_\infty /2.
\]
This is a contradiction, therefore IR condition and beta-min condition  $(\mathop{min}_{j \in S_0}| \beta^{0}_{j} | > \lambda \| \Sigma^{-1}_{11} \tau_S \|_\infty /2 )$ together implies variable selection ${S}^{\star} = S_0$ .
\qed
\end{pf}

\begin{theorem}[IR is necessary condition for variable selection]
Let ${S}^{\star}$ be the variables selected by the Lasso in (\ref{noiselessLasso}). Suppose ${S}^{\star} \subset S$ then the irrepresentable condition is met for the true active set $S$.
\end{theorem}
\begin{pf}
Given  ${S}^{\star} \subset S$ , it implies  $\beta_{2}  = 0( {\beta^{\star}_{S^c}}  = 0)$.
We get the following equation using the KKT condition.
\begin{align*}
	\left[ \begin{array}{cc}
	\Sigma_{11} & \Sigma_{12}\\
	\Sigma_{21} &\Sigma_{22})	
	\end{array} \right] 	 \left( \begin{array}{c} \beta_1- \beta^{0}_{1} \\ 0 \end{array} \right) & = \left( \begin{array}{c} \tau_1 \\ \tau_2 \end{array} \right) 
\end{align*}
After simplification we get
\begin{align*}
	2\Sigma_{11}(\beta_1-\beta^{0}_{1})  &= -\lambda \tau_1\\
	2\Sigma_{21}(\beta_1-\beta^{0}_{1}) &= -\lambda \tau_2.
\end{align*}
Multiplying the first equality with $\Sigma^{-1}_{11}$ we get
\begin{align*}
	2(\beta_1-\beta^{0}_{1}) &= -\lambda \Sigma^{-1}_{11} \tau_1.
\end{align*}
By substituting $ 2(\beta_1-\beta^{0}_{1})$ in equality 2 we get the following.
\begin{align*}
	\Sigma_{21}(\lambda \Sigma^{-1}_{11} \tau_1)  &= -\lambda \tau_2 \\
	 \| \Sigma_{21} \Sigma^{-1}_{11} \tau_1\|_{\infty} &=  \| {\beta}^{\star}_{S^c}\|_{\infty} \leq 1\\
	  \| \Sigma_{21} \Sigma^{-1}_{11} \tau_S\|_{\infty} & \leq 1
\end{align*}
\qed
\end{pf}

\subsubsection{The irrepresentable condition implies the compatibility
condition} \label{sec:IRRE}
We have shown that the irrepresentable condition implies variable selection, now we show that it is more restrictive than the compatibility condition. The following result (and proof) is from \cite{Sara3}.
\begin{theorem}
For a fixed set $S$, uniform $\theta$-IR condition implies compatibility condition.% with constant $\phi_{comp}(L,S) \geq 0$. 
\end{theorem}
\begin{pf}
It is given that uniform $\theta$-IR condition is satisfied by the design matrix for the set $S$, hence it is implicit that $\Lambda_{min}(\Sigma_{11}) > 0$. Now Let us suppose that compatibility condition does not hold for the set $S$, that is $\phi_{comp}^{2}(L,S) = 0$.\\
\\
The compatibility condition can also be defined as:
\begin{align*}
	\phi_{comp}^{2}(L,S) &= min \left\lbrace \frac{s}{n} \| \textbf{X} \Delta \|_{2}^{2} : \| \Delta_{S}\|_1 = 1, \; \|  \Delta_{S^c}\|_1 \leq L \right\rbrace,
\end{align*}
where  $\frac{1}{n} \| \textbf{X} \Delta \|_{2}^{2}  = \Delta' \Sigma \Delta$.
Suppose $\tilde{\Delta}$ solves the above minimization problem, that is
\begin{align}
	 \tilde{\Delta} = \mathop{arg min}_{\Delta \in \mathbb{R}^p} \left\lbrace \Delta' \Sigma \Delta : \| \Delta_{S}\|_1 = 1, \; \|  \Delta_{S^c}\|_1 \leq L \right\rbrace. \label{eq:minProblem}
\end{align}
We assume that the first $s$ variables are the active variables, and we partition the $\Sigma $ and $\tilde{\Delta} = (\tilde{\Delta}_1, \tilde{\Delta}_2)' $ accordingly. 
\begin{align*}
   \tilde{\Delta}' \Sigma \tilde{\Delta} & = (\tilde{\Delta}^{'}_{1}, \tilde{\Delta}^{'}_{2}) \left[ \begin{array}{cc}
	\Sigma_{11} & \Sigma_{12}\\
	\Sigma_{21} &\Sigma_{22})	
	\end{array} \right] 	 ( \begin{array}{c} \tilde{\Delta}_1 \\ \tilde{\Delta}^{'}_{2} \end{array}) 
\end{align*}
Under assumption of $\phi_{comp}^{2}(L,S) = 0$, the following equality holds.
\begin{align} \label{eq:cc}
	\tilde{\Delta} ' \Sigma \tilde{\Delta} &= 0 \\
	\tilde{\Delta}^{'}_{1} \Sigma_{11} \tilde{\Delta}_1+  \tilde{\Delta}^{'}_{1} \Sigma_{12} \tilde{\Delta}_2 + \tilde{\Delta}^{'}_{2} \Sigma_{21} \tilde{\Delta}_1 + \tilde{\Delta}^{'}_{2} \Sigma_{22} \tilde{\Delta}_2 &= 0 
\end{align}
We introduce a Lagrange multiplier $\lambda \in \mathbb{R}$ for the equality constraint $ \| \Delta_1 \|_1 = 1$ in Eq. (\ref{eq:minProblem}). Then by the KKT conditions, there exists a vector $\tau_1$ , such that $\| \tau_1 \|_\infty \leq 1$,  $\tau^{'}_{1} \Delta_1 =  \| \Delta \|_1$ and 
\begin{align} \label{eq:kkt1}
	\Sigma_{11} \tilde{\Delta}_1 + \Sigma_{12} \tilde{\Delta}_2 &= -\lambda \tau_1	.
\end{align}
By multiply  $(\tilde{\Delta}_1)' $ we obtain 
\begin{align}
	(\tilde{\Delta}_1)' \Sigma_{11} \tilde{\Delta}_1 +(\tilde{\Delta}_1)' \Sigma_{12} \tilde{\Delta}_2 &= -\lambda \| \Delta \|_1 \nonumber  \\
	(\tilde{\Delta}_1)' \Sigma_{11} \tilde{\Delta}_1 +(\tilde{\Delta}_1)' \Sigma_{12} \tilde{\Delta}_2 &= -\lambda.
\end{align}
By multiplying $(\tau_1)'\Sigma^{-1}_{11}$ in Eq. (\ref{eq:kkt1}) we get
\begin{align*}
	(\tau_1)' \tilde{\Delta}_1 +(\tau_1)'\Sigma^{-1}_{11} \Sigma_{12} \tilde{\Delta}_2 &= - \lambda (\tau_1)'\Sigma^{-1}_{11}  \tau_1 \\
	\| \tilde{\Delta}_1 \|_1 +(\tau_1)'\Sigma^{-1}_{11} \Sigma_{12} \tilde{\Delta}_2 &= - \lambda (\tau_1)'\Sigma^{-1}_{11}  \tau_1.
\end{align*}
By substituting $\| \tilde{\Delta}_1 \|_1  = 1$ we get
\begin{align}
	1 &= - (\tau_1)'\Sigma^{-1}_{11} \Sigma_{12} \tilde{\Delta}_2 -\lambda (\tau_1)'\Sigma^{-1}_{11} \tau_1. \label{eq:intermid1}
\end{align}
We can simplify the term $- (\tau_1)'\Sigma^{-1}_{11} \Sigma_{12} \tilde{\Delta}_2 $ as follows
\begin{align*}
- (\tau_1)'\Sigma^{-1}_{11} \Sigma_{12} \tilde{\Delta}_2  \leq |(\tau_1)'\Sigma^{-1}_{11} \Sigma_{12} \tilde{\Delta}_2| & \leq \|(\tau_1)'\Sigma^{-1}_{11} \Sigma_{12}\|_{\infty} \|\Delta_2 \|_{1} \leq L \theta\\ 	
	\because  \|(\tau_1)'\Sigma^{-1}_{11} \Sigma_{12}\|_{\infty} \leq  \theta & \text{ and } \|\Delta_2 \|_{1} \leq L .
\end{align*}
By substituting value for $- (\tau_1)'\Sigma^{-1}_{11} \Sigma_{12} \tilde{\Delta}_2$ in Eq. (\ref{eq:intermid1}) we obtain
\begin{align*}
 1- L \theta & \leq - \lambda (\tau_1)'\Sigma^{-1}_{11} \tau_1 .
\end{align*}
Let us assume that $\theta < \frac{1}{L}$. Then $1- L \theta > 0$ and $(\tau_1)'\Sigma^{-1}_{11} \tau_1 > 0$, therefore $\lambda < 0$, and we get the following inequality.
\begin{align}
	- \lambda \geq \frac{(1-L\theta) \Lambda_{min}(\Sigma_{11})}{s} \label{eq:lambda}
\end{align}
 Here, we used the inequality $(\tau_1)'\Sigma^{-1}_{11} \tau_1 \leq \frac{s}{\Lambda_{min}(\Sigma_{11})}$, which is due to \\
 $
  \frac{(\tau_1)'\Sigma^{-1}_{11} \tau_1 }{\|\tau \|^{2}_{2}} \leq \frac{1}{\Lambda_{min}(\Sigma_{11})}
 $.
By multiplying $( \tilde{\Delta}_2)' \Sigma_{21} \Sigma^{-1}_{11}$ in Eq. (\ref{eq:kkt1}) we get
\begin{align*}
	 (\tilde{\Delta}_2)' \Sigma_{21} \tilde{\Delta}_1+ ( \tilde{\Delta}_2)' \Sigma_{21} \Sigma^{-1}_{11} \Sigma_{12} \tilde{\Delta}_2 &= -\lambda ( \tilde{\Delta}_2)' \Sigma_{21} \Sigma^{-1}_{11} \tau_1.	
\end{align*}
Here, we consider projecting $(X_{S^{c}} \Delta_2)$ to the space spanned by $X_S$. The projected vector is $X_S \Sigma^{-1}_{11} \Sigma_{12} \tilde{\Delta}_2$. $\ell_2$ norm of the projected vector is 
\[
(X_S \Sigma^{-1}_{11} \Sigma_{12} \tilde{\Delta}_2)' X_S \Sigma^{-1}_{11} \Sigma_{12} \tilde{\Delta}_2 = ( \tilde{\Delta}_2)' \Sigma_{21} \Sigma^{-1}_{11} \Sigma_{12} \tilde{\Delta}_2\]
and $\ell_2$ norm of $(X_{S^{c}} \Delta_2)$ can also be written as \[(\tilde{\Delta}_2)' \Sigma_{22} \tilde{\Delta}_2. \]
Therefore $( \tilde{\Delta}_2)' \Sigma_{21} \Sigma^{-1}_{11} \Sigma_{12} \tilde{\Delta}_2 = (\tilde{\Delta}_2)' \Sigma_{22} \tilde{\Delta}_2$.
Now we use  the following fact.
\begin{align*}
 - \lambda ( \tilde{\Delta}_2)' \Sigma_{21} \Sigma^{-1}_{11} \tau_1 &= |\lambda| ( \tilde{\Delta}_2)' \Sigma_{21} \Sigma^{-1}_{11} \tau_1 \quad (\because \lambda < 0) \\
 & \leq |\lambda| |( \tilde{\Delta}_2)' \Sigma_{21} \Sigma^{-1}_{11} \tau_1 | \\
 & \leq |\lambda| L \theta \quad (\because |( \tilde{\Delta}_2)' \Sigma_{21} \Sigma^{-1}_{11} \tau_1 |  \leq L\theta) \\
  \lambda ( \tilde{\Delta}_2)' \Sigma_{21} \Sigma^{-1}_{11} \tau_1 & \leq - \lambda L \theta \\
	  - \lambda ( \tilde{\Delta}_2)' \Sigma_{21} \Sigma^{-1}_{11} \tau_1 & \geq \lambda L \theta
\end{align*}
Substituting value for $- \lambda ( \tilde{\Delta}_2)' \Sigma_{21} \Sigma^{-1}_{11} \tau_1$, we obtain
\begin{align} \label{eq:d22}
  (\tilde{\Delta}_2)' \Sigma_{21} \tilde{\Delta}_1 + (\tilde{\Delta}_2)' \Sigma_{22} \tilde{\Delta}_2 \geq \lambda L \theta.
\end{align}
Finally after substituting from Eq. (\ref{eq:kkt1}) and Eq. (\ref{eq:d22})  values to Eq. (\ref{eq:cc}), we get the following inequality.
\begin{align*}
	\tilde{\Delta} ' \Sigma \tilde{\Delta} = -\lambda + \lambda L \theta & \geq 0\\
	\therefore \tilde{\Delta} ' \Sigma \tilde{\Delta} & \geq 0
\end{align*}
Which contradicts our assumption that $\phi_{comp}(L,S) = 0$, and continuing with the above inequality, we get the bound for $\phi_{comp}(L,S)$ as follows.
\begin{align*}
	\phi_{comp}(L,S) &= s \tilde{\Delta} ' \Sigma \tilde{\Delta}\\
	\phi_{comp}(L,S) &\geq -\lambda s (1- L \theta)
\end{align*}
Substituting value of $ -\lambda$ from Eq. (\ref{eq:lambda}), we get the following bound.
\begin{align*}
	\phi_{comp}(L,S) \geq (1-L \theta)^2 \Lambda^{2}_{min}(\Sigma_{11}) 
\end{align*}
\qed
\end{pf}

\section{Conclusion}
We discussed various conditions required for the Lasso to hold oracle-inequalities. The compatibility condition is the weakest among others and we used it for deriving oracle results for prediction (fast rates) and estimation ($\ell_1$-norm). The oracle results for the slow rates for prediction does not assume any condition on the design matrix. We derived oracle results for $\ell_2$ estimation error using a slightly stronger version of the restricted eigenvalue condition. We discussed that the oracle result for variable selection requires irrepresentable condition and beta-min conditions. We illustrated various design conditions using simple examples. We also discussed that the irrepresentable condition implies compatibility condition. For further details on how various conditions for Lasso oracle results relate to each other, we refer to \cite{Sara3}.

\section*{References}
  \bibliographystyle{elsarticle-num} 
  \bibliography{niharika}

\end{document}